\documentclass[12pt,amssymb,verbatim]{amsart}

\usepackage[all]{xy}

\def\cal{\mathcal}

\numberwithin{equation}{section}

\theoremstyle{plain}
\newtheorem{theorem}{Theorem}

\newtheorem{proposition}[theorem]{Proposition}
\newtheorem{lemma}[theorem]{Lemma}

\newtheorem{corollary}[theorem]{Corollary}
\newtheorem{definition}[theorem]{Definition}

\theoremstyle{definition}

\theoremstyle{remark}
\newtheorem{remark}[theorem]{Remark}

 \def\today{\ifcase\month\or
  January\or February\or March\or April\or May\or June\or
  July\or August\or September\or October\or November\or December\fi
  \space\number\day, \number\year}

\begin{document}
\title[Outer forms of type $A_2$]
{Outer forms of type $A_2$ with infinite genus}
\author[Tikhonov]{Sergey V. Tikhonov}


\address{
Belarusian State University, Nezavisimosti Ave., 4,
220030, Minsk, Belarus} \email{tikhonovsv@bsu.by }

\def\cA{{\cal A}}
\def\cB{{\cal B}}
\def\cC{{\cal C}}
\def\cD{{\cal D}}
\def\cV{{\cal V}}
\def\cE{{\cal E}}
\def\cM{{\cal M}}
\def\cR{{\cal R}}
\def\M{{\rm M}}
\def\cS{{\cal S}}
\def\Symb{\textrm{Symb}}
\def\Gal{\textrm{Gal}}
\def\U{\textrm{U}}
\def\SU{\textrm{SU}}
\def\R{\textrm{R}}
\def\GL{\textrm{GL}}

\def\U{{\rm U}}
\def\SU{{\rm SU}}
\def\SL{{\rm SL}}
\def\op{{\rm op\;}}


\def\cor{\textrm{cor}}
\def\deg{\textrm{deg}}
\def\exp{\textrm{exp}}
\def\Gal{\textrm{Gal}}
\def\ram{\textrm{ram}}
\def\Spec{\textrm{Spec}}
\def\Proj{\textrm{Proj}}
\def\Perm{\textrm{Perm}}
\def\coker{\textrm{coker\,}}
\def\Hom{\textrm{Hom}}
\def\im{\textrm{im\,}}
\def\ind{\textrm{ind}}
\def\int{\textrm{int}}
\def\inv{\textrm{inv}}
\def\min{\textrm{min}}

\def\w{\widehat}

\begin{abstract}

Let $G$ be an absolutely almost simple algebraic group over a field $K$. The genus ${\bf gen}_K(G)$ of $G$ is the set of
$K$-isomorphism classes of $K$-forms $G'$ of $G$ that have the same $K$-isomorphism classes of maximal $K$-tori as $G$.
We construct an example of outer forms of type $A_2$ with infinite genus.

\end{abstract}

\maketitle


\def\dd{{\partial}}

\def\into{{\hookrightarrow}}

\def\emptyset{{\varnothing}}

\def\alp{{\alpha}}  \def\bet{{\beta}} \def\gam{{\gamma}}
 \def\del{{\delta}}
\def\eps{{\varepsilon}}
\def\kap{{\kappa}}                   \def\Chi{\text{X}}
\def\lam{{\lambda}}
 \def\sig{{\sigma}}  \def\vphi{{\varphi}} \def\om{{\omega}}
\def\Gam{{\Gamma}}  \def\Del{{\Delta}}  \def\Sig{{\Sigma}}
\def\ups{{\upsilon}}


\def\A{{\mathbb A}}
\def\F{{\mathbb F}}
\def\Q{{\mathbb{Q}}}
\def\CC{{\mathbb{C}}}
\def\PP{{\mathbb P}}
\def\R{{\mathbb R}}
\def\Z{{\mathbb Z}}
\def\X{{\mathbb X}}

\def\Gm{{{\Bbb G}_m}}
\def\Gmk{{{\Bbb G}_{m,k}}}
\def\GmL{{\Bbb G_{{\rm m},L}}}
\def\Ga{{{\Bbb G}_a}}

\def\Fb{{\overline F}}
\def\Hb{{\overline H}}
\def\Kb{{\overline K}}
\def\Lb{{\overline L}}
\def\Yb{{\overline Y}}
\def\Xb{{\overline X}}
\def\Tb{{\overline T}}
\def\Bb{{\overline B}}
\def\Gb{{\overline G}}
\def\Vb{{\overline V}}

\def\kb{{\bar k}}
\def\xb{{\bar x}}

\def\Th{{\hat T}}
\def\Bh{{\hat B}}
\def\Gh{{\hat G}}

\def\Xt{{\tilde X}}
\def\Gt{{\tilde G}}

\def\gg{{\mathfrak g}}
\def\gm{{\mathfrak m}}
\def\gp{{\mathfrak p}}
\def\gq{{\mathfrak q}}

\def\textrm#1{\text{\rm #1}}

\def\res{\textrm{res}}
\def\cor{\textrm{cor}}
\def\R{\textrm{R}}

\def\tors{_{\textrm{tors}}}      \def\tor{^{\textrm{tor}}}
\def\red{^{\textrm{red}}}         \def\nt{^{\textrm{ssu}}}
\def\sc{^{\textrm{sc}}}
\def\sss{^{\textrm{ss}}}          \def\uu{^{\textrm{u}}}
\def\ad{^{\textrm{ad}}}           \def\mm{^{\textrm{m}}}
\def\tm{^\times}                  \def\mult{^{\textrm{mult}}}
\def\tt{^{\textrm{t}}}
\def\uss{^{\textrm{ssu}}}         \def\ssu{^{\textrm{ssu}}}
\def\cf{^{\textrm{cf}}}
\def\ab{_{\textrm{ab}}}

\def\et{_{\textrm{\'et}}}
\def\nr{_{\textrm{nr}}}

\def\SB{\textrm{SB}}

\def\<{\langle}
\def\>{\rangle}

\def\til{\;\widetilde{}\;}


\font\cyr=wncyr10 scaled \magstep1%
\def\Bcyr{\text{\cyr B}}
\def\Sh{\text{\cyr Sh}}
\def\Ch{\text{\cyr Ch}}

\def\lpsi{{{}_\psi}}
\def\bks{{\backslash}}

\def\Br{\textrm{Br}}
\def\Pic{\textrm{Pic}}
\def\Bt{{{}_2\textrm{Br}}}
\def\Bn{{{}_n\textrm{Br}}}
\def\Grass{\textrm{Grass}}



Let $K$ be a field and $K^{sep}$ its separable closure.
Two absolutely almost simple algebraic $K$-groups $G_1$ and $G_2$ are said to have the same $K$-isomorphism classes of
maximal $K$-tori if every maximal $K$-torus of $G_1$ is $K$-isomorphic to some maximal $K$-torus of $G_2$,
and vice versa. An algebraic $K$-group $G'$ is called a $K$-form of an algebraic $K$-group $G$ if $G$ and $G'$ are
isomorphic over  $K^{sep}$.

\begin{definition} \cite[Def. 6.1]{ChRaRa15}
Let $G$ be an absolutely almost simple algebraic group over a field $K$. The genus ${\bf gen}_K(G)$ of $G$ is the set of
$K$-isomorphism classes of $K$-forms $G'$ of $G$ that have the same $K$-isomorphism classes of maximal $K$-tori as $G$.
\end{definition}

The genus is trivial in some special cases and it is conjectured to be finite whenever the field $K$ is finitely generated of "good" characteristic (see details in \cite[\S 8]{RaRa20}).

In a similar way one can define the genus of a division algebra.

\begin{definition}
The genus ${\bf gen}(\cD)$ of a finite-dimensional central division algebra $\cD$ over a field $K$ is defined as the set of classes $[\cD']\in \Br(K)$,
where $\cD'$ is a central division $K$-algebra having the same maximal subfields as $\cD$.
\end{definition}

If $\cD$ is a finite-dimensional central division $K$-algebra, then it is well-known that any maximal $K$-torus of the corresponding algebraic group $G=\SL_{1,\cD}$
is of the form $\R_{E/K}(\mathbb{G}_m) \cap G$ (where $\R_{E/K}(\mathbb{G}_m)$ is the Weil restriction of the 1-dimensional split torus $\mathbb{G}_m$) for some maximal separable subfield $E$ of $\cD$.
Thus the results on genus of division algebras from \cite{Me14} and \cite{Ti16}
rephrased in the language of algebraic groups say that for any prime $p$, there exist fields (with infinite transcendence degree over the prime subfield)
over which there are  inner forms of type $A_{p-1}$ with infinite genus.
An example of  groups of type $G_2$ with infinite genus is obtained in \cite[Rem. 3.6(b)]{BeGiLe16}.
In the present paper, we construct such an example for outer forms of type $A_2$.

Let $F/K$ be a quadratic separable field
extension and $\sigma$ the non-trivial $K$-automorphism of $F$. An involution on an $F$-algebra $\cR$ is
called an $F/K$-involution if its restriction to $F$ is $\sigma$.
An isomorphism of $F$-algebras with involution $f: (\cR,\tau) \longrightarrow (\cR',\tau')$ is an $F$-algebra isomorphism
$f: \cR \longrightarrow \cR'$ such that $\tau' \circ f = f \circ \tau$.
Let also $E/F$ be a field extension such that $E$ has an automorphism $\psi$ of order 2
such that $\psi | F = \sigma$ (i.e., $E$ has an $F/K$-involution). Then $\tau_E$ denotes the involution on $\cR \otimes_F E $ defined by the formula
${\tau_E} (r\otimes e) := \tau(r) \otimes \psi(e)$, where $r\in \cR$, $e \in E$. In particular, if $L/F$ is a field extension
linearly disjoint to $E$ over $F$ with an
automorphism $\phi$ of order two extending $\sigma$, then $\phi_E$ is the automorphism of order two of the field $EL=L\otimes_F E$ which
extends $\phi$.

Let $\cA$ be a central division $F$-algebra of degree $n$ with an $F/K$-involution $\tau$.
Let $L/F$ be a separable field extension of degree $n$, and let $\phi: L \longrightarrow L $ be an automorphism of order two such that
$\phi| F = \sigma$. An embedding of algebras with involution $(L,\phi) \hookrightarrow (\cA,\tau)$ is by definition an injective
$F$-homomorphism $f:L\longrightarrow \cA$ such that $\tau \circ f =f \circ \phi$.
It is known that embeddings of maximal tori into the special unitary group $\SU(\cA,\tau)$ can be described in terms of
embeddings of fields with involution into the central simple algebra with involution $(\cA,\tau)$ (\cite[Prop. 2.3]{PrRa10}).

In this paper, we construct a field
$E$ and a subfield $T\subset E$ such that $[E:T]=2$ and
there is an infinite set of (pairwise non-isomorphic) division $E$-algebras of degree 3 with $E/T$-involution such that
a field extension $L/E$ of degree 3 can be embedded as an algebra with involution into one algebra of this set if and only
if it can be embedded as an algebra with involution into all other algebras from this set. Passing to the corresponding special unitary groups, we obtain an example of outer
forms of type $A_2$ with infinite genus.


Below we use the following notation:
$Alg_3(F/K)$ is the set of isomorphism classes of central division $F$-algebras of degree $3$ with $F/K$ involution;
$Ext_3(F/K)$ is the set of isomorphism classes of field extensions of $F$ of degree $3$  with $F/K$-involution.
The $3$-torsion of the Brauer group $\Br(F)$ is denoted by $_3\Br(F)$.
For a field extension $E/F$ and a central simple $F$-algebra $\cA$, $\cA_E$ denotes the tensor product
$\cA\otimes_F E$ and $\res_{E/F} : \Br(F)\longrightarrow \Br(E)$ denotes the restriction homomorphism. The restriction of $\res_{E/F}$ to the subgroup
$_3\Br(F)$ will also be denoted by $\res_{E/F}$. For a central simple $F$-algebra $\cA$, $\cA^{op}$ denotes the opposite algebra and
$\cA^m$ denotes $\cA\otimes_F \dots \otimes_F \cA$ ($m$ times).
For a quadratic form $q$ over $K$ and a field extension $E/K$,  $q_E$ denotes the quadratic form obtained by
extension of scalars from $K$ to $E$.
Recall that a field extension $E/F$ is called regular if $E/F$ is separable and $F$ is algebraically closed in $E$.


We start with the following


\begin{lemma} \label{l:one_algebra_cyclic_extension}
Let $n$ be a positive integer, $F$ a field of characteristic not dividing $2n$, $F/K$ a quadratic field
extension, $\sigma$ the non-trivial $K$-automorphism of $F$, $\cA$ a central simple
$F$-algebra of degree $n$, and
$L/F$ a cyclic field extension of degree $n$. Then there exists a regular field extension $M/F$ and
a subfield $T\subset M$ such that $[M:T]=2$ and

(1) $M = T F$ and the non-trivial $T$-automorphism of $M$ extends $\sigma$;


(2) the composite $M L$ splits $\cA_M$;

(3) the homomorphism $\res_{M/F} : \Br(F) \longrightarrow \Br(M)$ is injective.
\end{lemma}

\noindent {\it {Proof}}. Let  $F(x)$ be a purely transcendental extension of $F$ of transcendence degree 1.
Let also $\phi$ be a generator of the Galois group $\Gal(L(x)/F(x))$ and
$$
\cC := \cA_{F(x)}^{op}\otimes_{F(x)} (L(x)/F(x),\phi,x),
$$
where $(L(x)/F(x),\phi,x)$ is a cyclic $F(x)$-algebra of degree $n$.
Let also $E$ be the function field of the Severi-Brauer variety of $\cC$.
Note that the kernel of the restriction homomorphism
$\res_{E/F(x)}: \Br(F(x)) \longrightarrow \Br(E)$ is generated by
$[\cC]$ (see, e.g., \cite[Cor. 13.16]{Sa99}).

Let $\cB$ be a central simple  $F$-algebra of exponent bigger than 1. Assume that $\cB$ is split by $E$, then $[\cB_{F(x)}]=[\cC^i]$ for some $1\le i \le n$.
If $i<n$, then the $F(x)$-algebra $\cC^i$ ramifies at the discrete valuation (trivial on $F$) of $F(x)$ defined by the polynomial $x$, but $\cB_{F(x)}$ is unramified
at this valuation, hence $[\cB_{F(x)}]\ne [\cC^i]$. Since the exponent of $\cB_{F(x)}$ is bigger than 1, then $[\cB_{F(x)}]\ne [\cC^n]=[F(x)]$. Thus $\cB_{F(x)}$ is not split by $E$, i.e.,
the homomorphism $\res_{E/F} : {\Br(F)} \longrightarrow {\Br(E)}$ is injective.

Since $E$ splits $\cC$, then
$$
[\cA_{E}] = [(L(x)/F(x),\phi,x)_{E}]= [(E L/E,\phi',x)],
$$
where $\phi'$ is the generator of $\Gal(E L/E)$. Thus $E L$ splits $\cA_{E}$.

Note that $E/F$ is a regular extension of $F(x)$.
For the following construction
of the transfer of a regular field extension, we refer to \cite[p. 220]{Ro64}).

Let $\sigma$ also denotes the $K(x)$-automorphism of $F(x)$  extending the automorphism
$\sigma$ of $F$.
The automorphism $\sigma$ of $F(x)$ can be extended to an isomorphism (which we also denote by $\sigma$) of $E$ and another regular extension of $F(x)$
denoted by $E_{\sigma}$. Thus the following diagram commutes:

$$
\xymatrix{F \ar@{^{(}->}[r] \ar[d]^-{\sigma} & F(x) \ar@{^{(}->}[r] \ar[d]^-{\sigma}&
E \ar[d]^-{\sigma}\\
F \ar@{^{(}->}[r] & F(x) \ar@{^{(}->}[r] &  E_{\sigma} \rlap{\ .}\\
}
\eqno{(1)}
$$

Let $M=E E_{\sigma}$ be the free composite over $F$ of $E$ and
$E_{\sigma}$. This free composite is $F$-isomorphic to the function field of the Severi-Brauer varity of
the $E_{\sigma}(y)$-algebra $\cA_{E_{\sigma}(y)}^{op}\otimes_{E_{\sigma}(y)} (E_{\sigma}L(y)/E_{\sigma}(y),\psi',y)$, where $y$ is transcendental over $E_{\sigma}$
(we replace $x$ by $y$ since the composite is free) and $\psi'$ is the generator of the Galois group $\Gal(E_{\sigma}L(y)/E_{\sigma}(y))$.
The field $M$ is a regular extension of $F$.
The isomorphisms $\sigma : E \longrightarrow E_{\sigma}$ and $\sigma^{-1} : E_{\sigma} \longrightarrow E$ have a unique extension
to an automorphism $\bar \sigma$ of $M$ of order two.
Let $T:=T_{F/K}(E)$ be the transfer of $E$ with respect to the ground
field descent $F\supset K$, i.e., $T$ is the subfield of $M$ of elements
fixed under the action of $\bar \sigma$.
Note that the composite $TF$ coincides with $M$, $[M:T]=2$, and $\bar \sigma$ extends $\sigma$.

The algebra $\cA_M$ is split by $ML$ since $\cA_E$ is split by $EL$.


Finally, the diagram (1) induces the following commutative diagram for the corresponding Brauer groups:
$$
\xymatrix{\Br(F) \ar[r]^-{\res} \ar[d]^-{\cong} & \Br(F(x)) \ar[r]^-{\res} \ar[d]^-{\cong}&
\Br(E) \ar[d]^-{\cong}\\
\Br(F) \ar[r]^-{\res} & \Br(F(x)) \ar[r]^-{\res} &  \Br(E_{\sigma}) \ar[r]^-{\res} &  \Br(M) \rlap{\ .}\\
}
$$
The injectivity of $\res_{E/F}$ implies the injectivity of $\res_{E_{\sigma}/F}$. Moreover, $\res_{M/E_{\sigma}}$ is injective
by the same arguments as for $\res_{E/F}$, we just replace the ground field $F$ by $E_{\sigma}$. Hence the homomorphism $\res_{M/F}$ is also injective.

\qed


We also need the following


\begin{lemma} \label{l:fields}
Let $F$ be a field of characteristic $\ne 2, 3$; $F/K$ a quadratic field
extension, $\sigma$ the non-trivial $K$-automorphism of $F$, and
$L/F$ a field extension of degree $3$ with an automorphism $\phi$ of order two
such that $\phi |F = \sigma$. Then there exists a field extension $F(L)/F$  and
a subfield $K(L)\subset F(L)$ such that $[F(L):K(L)]=2$ and

(1) $F(L) = K(L) F$ and the non-trivial $K(L)$-automorphism of $F(L)$, denoted by $\sigma_{F(L)}$, extends $\sigma$;

(2)  $[F(L): F] \le 2$; 

(3)  the composite $F(L) L$ is a cyclic extension of $F(L)$ of degree 3;

(4)  the homomorphism $\res_{F(L)/F} : {_3\Br(F)} \longrightarrow {_3\Br(F(L))}$ is injective. 
\end{lemma}


\noindent {\it {Proof}}.
If the extension $L/F$ is cyclic, then one can take $F(L) := F$, $K(L) := K$, and $\sigma_{F(L)} := \sigma$.

Assume that the extension $L/F$ is not cyclic. Let $N$ be the normal closure of the extension $L^{\phi}/K$, where $L^{\phi}\subset L$ is the subfield of elements fixed by $\phi$.
Then $F \not \subset N$ and $NF$ is the normal closure of the extension $L/F$.
Thus we have the following diagram of field extensions:

$$
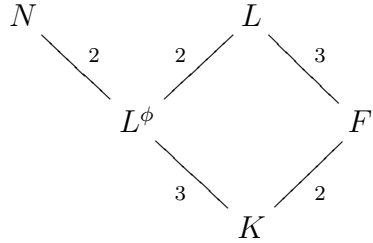
\xymatrix{
       N &   & L &\\
& \ar@{-}[lu]_2 L^{\phi}\ar@{-}[ru]^2 & & F \ar@{-}[lu]_3 \\
&& \ar@{-}[lu]^3  K \ar@{-}[ru]_2 & &  \\
}
$$

Let $H$ be the Sylow 3-subgroup of the Galois group $\Gal(N/K)$. Then $N^H$, the fixed field of $H$, is an extension
of $K$ of degree 2 and $N/N^H$ is a cyclic extension of degree 3.  Hence $NF$ is a cyclic extension of $N^HF$ of degree 3
and $[N^HF : F] = 2$. Let $K(L) := N^H$ and $F(L) := K(L)F$.
Since $F\not\subset N$, then $[F(L):K(L)]=2$ and the field $F(L)$
has a $K(L)$-automorphism of order two extending $\sigma$. Note that $F(L)L=NF$, hence $F(L) L/F(L)$ is a cyclic extension of degree 3.
Finally, since $[F(L):F] = 2$, then the homomorphism $\res_{F(L)/F} : {_3\Br(F)} \longrightarrow {_3\Br(F(L))}$ is injective.
\qed


\begin{remark}
In the notations of Lemma \ref{l:fields}, for any field extension $L'$ of $F$ of degree 3, $F(L)$ and $L'$
are linearly disjoint over $F$. Moreover, if $L'$ has an automorphism $\phi'$ of order 2 extending $\sigma$,
then the composite $F(L)L'$ has the automorphism $\phi'_{F(L)}$ which extends the automorphisms $\phi'$ and $\sigma_{F(L)}$.
\end{remark}


\begin{lemma} \label{l:one_algebra_any_extension}
Let $F$ be a field of characteristic $\ne 2, 3$; $F/K$ a quadratic field
extension, $\sigma$ the non-trivial $K$-automorphism of $F$, $\cA$ a central simple
$F$-algebra of degree $3$ with an $F/K$-involution $\tau$, and
$L/F$ a field extension of degree $3$ with an automorphism $\phi$ of order two
such that $\phi |F = \sigma$. Then there exists a field extension $F_{(K,L,\cA)}/F$
and a subfield $K_{(L,\cA)}\subset F_{(K,L,\cA)}$
such that $[F_{(K,L,\cA)}:K_{(L,\cA)}]=2$ and


(1) $F_{(K,L,\cA)} = K_{(L,\cA)} F$ and the non-trivial $K_{(L,\cA)}$-automorphism of $F_{(K,L,\cA)}$, denoted by $\sigma_{F_{(K,L,\cA)}}$, extends $\sigma$;

(2) the homomorphism $\res_{F_{(K,L,\cA)}/F} : {_3\Br(F)} \longrightarrow {_3\Br(F_{(K,L,\cA)})}$ is injective;

(3) for any field extension $L'$ of $F$ of degree 3, $F_{(K,L,\cA)}$ and $L'$ are linearly disjoint over $F$;


(4) there is an embedding $(F_{(K,L,\cA)} L, \phi_{F_{(K,L,\cA)}}) \hookrightarrow
(\cA_{F_{(K,L,\cA)}}, \tau_{{F_{(K,L,\cA)}}} )$ of algebras with involution.
\end{lemma}


\noindent {\it {Proof}}.
Let $F(L)$, $K(L)$ and $\sigma_{F(L)}$ be as in Lemma \ref{l:fields}.
Let also $M$ and $T$ be fields obtained by applying Lemma \ref{l:one_algebra_cyclic_extension}
for the quadratic field extension $F(L)/K(L)$, the $F(L)$-algebra $\cA_{F(L)}$,
the cyclic field extension $F(L)L/F(L)$ of degree $3$.
Then by Lemmas \ref{l:one_algebra_cyclic_extension} and \ref{l:fields},
the homomorphism $\res_{M/F} : {_3\Br(F)} \longrightarrow {_3\Br(M)}$ is injective; for any field extension $L'$ of $F$ of degree 3,
$M$ and $L'$ are linearly disjoint over $F$ and the composite $M L$ splits $\cA_{M}$.
Thus there is an $M$-embedding $\varepsilon: ML \hookrightarrow \cA_M$ of $M$-algebras.

Note that $\cA_M$ has the $M/T$-involution $\tau_M$
which extends $\tau$ and $ML$ has the automorphism $\phi_M$ of order two extending $\phi$. By \cite[Proposition 3.1]{PrRa10},  there exists an $M/T$-involution $\delta$ on $\cA_M$ such that
$\varepsilon: (ML, \phi_M) \hookrightarrow (\cA_M, \delta)$ is an embedding of algebras with involution.

Let $\pi(\delta)$ and $\pi(\tau_M)$ be the 3-fold Pfister forms of involutions $\delta$ and $\tau_M$ respectively (see \cite[\S 19.B]{KnMeRoTi98}).
Let $T(\pi(\delta))$ and $T(\pi(\tau_M))$ be the function fields of $\pi(\delta)$ and $\pi(\tau_M)$ respectively.
Then the quadratic forms $\pi(\delta)_{T(\pi(\delta))}$ and $\pi(\tau_M)_{T(\pi(\tau_M))}$ are isotropic and hence hyperbolic since they are Pfister forms.

Let $K_{(L,\cA)}$ be the free composite over $T$ of the fields $T(\pi(\delta))$ and $T(\pi(\tau_M))$. Let also
$F_{(K,L,\cA)} := K_{(L,\cA)} F$. Since $F \not \subset K_{(L,\cA)}$, then $[F_{(K,L,\cA)}:K_{(L,\cA)}]=2$ and $F_{(K,L,\cA)}$ has a
$K_{(L,\cA)}$-automorphism $\sigma_{F_{(K,L,\cA)}}$ of order 2 extending $\sigma$.
Note that the algebraic closure of $F$ in $F_{(K,L,\cA)}$ is $F(L)$ and $[F(L):F]$ is either 1 or 2. Therefore, $F_{(K,L,\cA)}$ and $L'$ are linearly disjoint over $F$
for any field extension $L'$ of $F$ of degree 3.


The extensions $T(\pi(\delta))/T$ and $T(\pi(\tau_M))/T$ are composition of a purely transcendental extension with a quadratic extension,
hence the homomorphism $\res_{F_{(K,L,\cA)}/M} : {_3\Br(M)} \longrightarrow {_3\Br(F_{(K,L,\cA)})}$ is injective. Hence
$\res_{F_{(K,L,\cA)}/F} : {_3\Br(F)} \longrightarrow {_3\Br(F_{(K,L,\cA)})}$ is also injective.

The quadratic forms $\pi(\delta)_{K_{(L,\cA)}}$ and $\pi(\tau_M)_{K_{(L,\cA)}}$ are hyperbolic. Then by \cite[Theorem 19.6]{KnMeRoTi98}),
the involutions $\delta_{F_{(K,L,\cA)}}$ and $\tau_{F_{(K,L,\cA)}}$ on $\cA_{F_{(K,L,\cA)}}$ are conjugate. This means that
there is an isomorphism $\xi:(\cA_{F_{(K,L,\cA)}}, \delta_{F_{(K,L,\cA)}}) \longrightarrow (\cA_{F_{(K,L,\cA)}}, \tau_{{F_{(K,L,\cA)}}} )$ of algebras with involution.

Moreover, the  embedding $\varepsilon: (ML, \phi_M) \hookrightarrow (\cA_M, \delta)$ of algebras with involution induces
an embedding $(F_{(K,L,\cA)} L, \phi_{F_{(K,L,\cA)}}) \hookrightarrow (\cA_{F_{(K,L,\cA)}}, \delta_{F_{(K,L,\cA)}})$ of algebras with involution.
Indeed, $F_{(K,L,\cA)} L = ML \otimes_M  F_{(K,L,\cA)}$. Let
$$
\varepsilon_{F_{(K,L,\cA)}} : ML \otimes_M  F_{(K,L,\cA)} \longrightarrow \cA_{F_{(K,L,\cA)}}
$$
be an $F_{(K,L,\cA)}$-embedding defined by the formula
$\varepsilon_{F_{(K,L,\cA)}}(m\otimes a) := \varepsilon(m)\otimes a$, where $m\in ML$, $a\in F_{(K,L,\cA)}$. Then
$$
\varepsilon_{F_{(K,L,\cA)}}( \phi_{F_{(K,L,\cA)}} (m\otimes a)) = \varepsilon_{F_{(K,L,\cA)}}( \phi_M(m)\otimes \sigma_{F_{(K,L,\cA)}}(a))=
\varepsilon(\phi_M(m) )\otimes \sigma_{F_{(K,L,\cA)}}(a) = $$
$$ \delta(\varepsilon(m) )\otimes \sigma_{F_{(K,L,\cA)}}(a) =
\delta_{F_{(K,L,\cA)}}(\varepsilon(m) \otimes a) = \delta_{F_{(K,L,\cA)}}(\varepsilon_{F_{(K,L,\cA)}} (m \otimes a)).
$$
Thus $\varepsilon_{F_{(K,L,\cA)}}$ is an embedding of algebras with involutions.
Then $\xi\circ \varepsilon_{F_{(K,L,\cA)}}$ is an embedding
$(F_{(K,L,\cA)}L,  \phi_{F_{(K,L,\cA)}}) \hookrightarrow (\cA_{F_{(K,L,\cA)}}, \tau_{{F_{(K,L,\cA)}}} )$ of algebras with involution.

\qed


The following construction of the field $F_{(K,S,A)}$ is an adaptation of the construction from \cite{Ti16} for algebras with involutions.
We give the details below for the reader's convenience.


\begin{proposition} \label{pr: fields}
Let $F$ be a field of characteristic $\ne 2, 3$; $F/K$ a quadratic field
extension, $\sigma$ the non-trivial $K$-automorphism of $F$, $A\subset Alg_3(F/K)$ and
$S\subset Ext_3(F/K)$. Then there exists a field extension $F_{(K,S,A)}/F$ and
a subfield $K_{(S,\cA)}\subset F_{(K,S,A)}$ such that $[F_{(K,S,A)}:K_{(S,A)}]=2$ and

(1) $F_{(K,S,A)} = K_{(S,A)} F$ and the non-trivial $K_{(S,A)}$-automorphism, denoted by $\sigma_{F_{(K,S,A)}}$, of $F_{(K,S,A)}$ extends $\sigma$;

(2) the homomorphism $\res_{F_{(K,S,A)}/F} : {_3\Br(F)} \longrightarrow {_3\Br(F_{(K,S,A)})}$ is injective;

(3) for any field extension $L'$ of $F$ of degree 3, $F_{(K,L,\cA)}$ and $L'$ are linearly disjoint over $F$;


(4) for any $\cA\in A$ with an $F/K$-involution $\tau$ and $L\in S$ with a $K$-automorphism $\phi$ of order 2 extending $\sigma$, there is an embedding $(F_{(K,S,A)} L, \phi_{F_{(K,S,A)}}) \hookrightarrow (\cA_{F_{(K,S,A)}}, \tau_{{F_{(K,S,A)}}} )$ of algebras with involution.

\end{proposition}


\noindent {\it {Proof}}. 
Let ${\mathcal{P}} :=\{(L,\cD) | L\in S \mbox{ and } \cD\in A\}$ be the set of pairs.
Let also $<$ be a well-ordering on $\mathcal{P}$ and let $t_0=(L_0,\cD_0)$ denote its least element.
Set $E_{t_0} := F_{(K, L_0,\cD_0)}$ and $T_0 := K_{(L_0,\cD_0)}$, where the fields $F_{(K, L_0,\cD_0)}$ and $K_{(L_0,\cD_0)}$ are constructed in
Lemma \ref{l:one_algebra_any_extension}. For $t=(L,\cD)\in {\mathcal{P}}$, set
$$
E^{<t} := \bigcup_{t' < t } E_{t'},  T^{<t} := \bigcup_{t' < t } T_{t'},  T_t : ={T^{<t}}_{(E^{<t} L,\cD_{E^{<t}})},
\mbox{ and } E_t : ={E^{<t}}_{(T^{<t}, E^{<t} L,\cD_{E^{<t}})},
$$
where the fields $E_t$ and $T_t$ are obtained by applying Lemma \ref{l:one_algebra_any_extension} to the quadratic field extension
$E^{<t}/T^{<t}$, the field extension $E^{<t} L/E^{<t}$ of degree 3, the automorphism $\phi_{E^{<t}}$ of $E^{<t} L$ extending the automorphism $\phi$ of $L$
and the $E^{<t}$-algebra $\cD_{E^{<t}}$.
We also define
$F_{(K,S,A)} := \bigcup_{t \in \mathcal{P} } E_t$ and $K_{(S,A)} := \bigcup_{t \in \mathcal{P} } T_t$.

By Lemma \ref{l:one_algebra_any_extension}, $E_t=T_t F$ and $[E_t : T_t]=2$ for any $t\in \mathcal{P}$. 
Then $F_{(K,S,A)} = K_{(S,A)} F$, $[ F_{(K,S,A)} : K_{(S,A)}] = 2$  and the non-trivial  $K_{(S,A)}$-automorphism of  $F_{(K,S,A)}$  extends $\sigma$.

By Lemma \ref{l:one_algebra_any_extension} and transfinite induction,
the homomorphism $\res_{F_{(K,S,A)}/F} : {_3\Br(F)} \longrightarrow {_3\Br(F_{(K,S,A)})}$ is injective and
for any field extension $L'$ of $F$ of degree 3, $F_{(K,S,A)}$ and $L'$ are linearly disjoint over $F$.

Finally, let $\cA\in A$ with an $F/K$-involution $\tau$, $L\in S$ with an automorphism $\phi$ of order 2 extending $\sigma$ and $t=(L,\cA)$.
By Lemma \ref{l:one_algebra_any_extension}, there is an embedding $(E_t L, \phi_{E_t}) \hookrightarrow (\cA_{E_t}, \tau_{E_t} )$ of algebras with involution.
Moreover, as in the proof of Lemma \ref{l:one_algebra_any_extension}, this embedding induces
the embedding
$$
(F_{(K,S,A)} L, \phi_{F_{(K,S,A)}}) \hookrightarrow  (\cA_{F_{(K,S,A)}}, \tau_{{F_{(K,S,A)}}} )
$$
of algebras with involution.
\qed


\begin{theorem}  \label{th: subset}
Let $F$ be a field of characteristic $\ne 2, 3$,
$F/K$ a quadratic field
extension, $\sigma$ the non-trivial $K$-automorphism of $F$, $A\subset Alg_3(F/K)$. Then there exists a field extension $F_A/F$ and
a subfield $K_A \subset F_A$ such that $[F_A:K_A]=2$ and

(1) $F_A = K_A F$ and the non-trivial $K_A$-automorphism, denoted by  $\sigma_{F_A}$, of $F_A$ extends $\sigma$;

(2) the homomorphism $\res_{F_A/F} : {_3\Br(F)} \longrightarrow {_3\Br(F_A)}$ is injective;

(3) for any central simple
$F$-algebra $\cB$ of degree $3$ with an $F/K$-involution $\theta$, the algebra $\cB_{F_A}$
has an $F_A/K_A$-involution $\theta_{F_A}$ extending $\theta$;

(4) if $L\in Ext_3(F_A/K_A)$ with a $K_A$-automorphism $\phi$ of order 2 extending $\sigma_{F_A}$,
then there is an embedding $(L, \phi) \hookrightarrow (\cA_{F_A}, \tau_{F_A} )$ of algebras with involution  for any $\cA\in A$ with an $F/K$-involution $\tau$.
\end{theorem}


\noindent {\it {Proof}}.
Let $F_0 := F$ and $K_0 := K$. We recursively define $F_i$ and $K_i$, $i\in {\mathbb{Z}}_ { > 0}$, to be the fields
${F_{i-1}}_{(K_{i-1}, Ext_3(F_{i-1}/K_{i-1}), \res_{F_{i-1}/F}(A))}$  and  ${K_{i-1}}_{(Ext_3(F_{i-1}/K_{i-1}), \res_{F_{i-1}/F}(A))}$
constructed by applying Proposition \ref{pr: fields} to the quadratic field extension $F_{i-1}/K_{i-1}$, the set
$\res_{F_{i-1}/F}(A)\subset Alg_3(F_{i-1}/K_{i-1})$ and
the set $Ext_3(F_{i-1}/K_{i-1})$.

Let $F_A := \bigcup_{i\ge 0} F_i$ and $K_A := \bigcup_{i\ge 0} K_i$.
Hence $F_A = K_A F$ and the non-trivial $K_A$-automorphism $\sigma_{F_A}$ of $F_A$ extends $\sigma$.
Therefore,  for any central simple $F$-algebra $\cB$ of degree $3$ with an $F/K$-involution $\theta$,
the $F_A/K_A$-involution $\theta_{F_A}$ extends $\theta$.

By induction and Proposition \ref{pr: fields},
$\res_{F_A/F} : {_3\Br(F)} \longrightarrow {_3\Br(F_A)}$ is injective.

Assume that $\cA \in A$ with an $F/K$-involution $\tau$ and $L \in Ext_3(F_A/K_A)$ with an automorphism $\phi$ of order two
extending  $\sigma_{F_A}$.
Then there exists $i\ge 0$ and   a field extension $L'$ of $F_i$ of degree 3 such that
$L= F_AL'$ and  $\phi_i :=  \phi_{|L'}$ is a $K_i$-automorphism of order two. This means that $L' \in Ext_3(F_i/K_i)$. By Proposition \ref{pr: fields},
there is an embedding
$$
(F_{i+1}L', {\phi_i}_{F_{i+1}}) \hookrightarrow (\cA_{{F_{i+1}}}, \tau_{F_{i+1}} )
$$
of algebras with involution.
As in the proof of Lemma \ref{l:one_algebra_any_extension}, this embedding can be extended to an embedding
$(L, \phi) \hookrightarrow (\cA_{F_A}, \tau_{F_A} )$ of algebras with involution.

\qed

As a corollary to Theorem \ref{th: subset}, we obtain the following

\begin{corollary} \label{th: infinite}
There exists a field $E$ and a subfield $T\subset E$ with $[E:T]=2$ such that there is an infinite set $B$
of pairwise non-isomorphic division $E$-algebras of degree 3 with $E/T$-involution and
such that for any field extension $L/E$ of degree 3 with an automorphism $\phi$ of order 2 extending the non-trivial $T$-automorphism of $E$,
there is an embedding $(L, \phi) \hookrightarrow (\cA, \tau)$ of algebras with involution for any $\cA \in B$ with an $E/T$-involution $\tau$.
\end{corollary}

\noindent {\it {Proof}}.
Let $\xi_3$ be a primitive 3th root of unity, $K=\mathbb{Q}(\xi_3)(x,y,z)$, the purely transcendental extension of the field $\mathbb{Q}(\xi_3)$ and
$F=K(\sqrt{2})$. Then for $i>0$, the symbol $F$-algebras $(z,\frac{x+\sqrt{2}y^i}{x-\sqrt{2}y^i})_3$ of degree 3 are pairwise non-isomorphic (since they
have different ramification) and have $F/K$-involutions (since the corestriction to $K$ of these algebras is trivial). Now we apply Theorem \ref{th: subset}
for the infinite set $A \subset Alg_3(F/K)$ consisting of these algebras and set $E := F_A$, $T := K_A$, $B := \res_{F_A/F}(A)$.

\qed


Rephrasing the previous corollary in the language of algebraic groups, we obtain the following

\begin{corollary} \label{th: infinite}
There exists a field $T$ such that there are  infinitely many (pairwise non-isomorphic)  outer forms of type $A_2$ over $T$ having the same
infinite genus.
\end{corollary}




\end{document}